\def\R{{\rm I}\!{\rm R}}

\parindent = 0cm



\centerline {\bf Lower bounds for pseudodifferential operators
with  a radial symbol. }

\bigskip

\centerline {Laurent Amour, Lisette Jager and Jean Nourrigat}

\medskip

\centerline {Universit\'e de Reims}

\vskip 1cm

\bigskip

\vbox {\sevenrm
ABSTRACT. In this paper we establish explicit lower bounds for pseudodifferential
operators with a radial symbol.
 The proofs use classical Weyl calculus techniques
 and some useful, if not celebrated, properties of the Laguerre polynomials.
}
\bigskip
 {\sevenrm
2010  Mathematical Subject Classification  Primary 35S05; Secondary 28C20,35R15,81S30.

  }

\bigskip 

{\bf 1. Introduction.}

\bigskip

If a function $F$ defined on $\R^{2d}$ is smooth and has bounded derivatives, 
 the Weyl calculus associates with it a pseudodifferential operator
 $Op_h^{Weyl}(F)$ which is bounded on   $L^2 (\R^d)$ and satisfies, for all
 $f$ and $g$ in  ${\cal S}(\R^d)$,
 $$  <  Op_h^{Weyl}(F)  f ,  g> =(2\pi h) ^{-d}
  \int _{\R^{2d}} F(Z)H_h(f , g , Z) dZ, \leqno (1.1)$$
where $H_h (f , g , \cdot )$ is the Wigner function

 $$H_h  (f , g , Z) = \int _{\R^d} e^{-{i\over h } t\cdot \zeta } f \left (
 z + {t\over 2} \right ) \overline {g \left (
 z - {t\over 2} \right )} dt \hskip 2cm Z = (z , \zeta) \in \R^{2d}. \leqno (1.2)$$
For this form of the definition, see [U], [L]  or  [C-R], Chapter II,
 Proposition 14.

 \bigskip

The different variants of  G\aa rding's inequality prove that, if
$F\geq 0$, the operator $ Op_h^{Weyl}(F)$ is roughly $\geq 0$. More precisely,
according to the classical  G\aa rding's inequality (see [HO] or [L]),
the non negativity of $F$ implies the existence of a positive constant $C$,
independent of $h$, 
such that, for all sufficiently small $h$ and for all $f$ in  ${\cal S } (\R^d)$:
 $$ < Op_h^{Weyl}(F)f ,f>\ \geq - C h \Vert f \Vert _{L^2(\R^d)}^2 .\leqno (1.3)$$
See [L-N] for other similar results.
This inequality holds for systems of operators, whereas the more precise
  Fefferman-Phong inequality [F-P] is valid only for scalar operators. 
 Fefferman-Phong's inequality states that, under the same hypotheses
as  G\aa rding's inequality, one has,  for all $h$ in $(0, 1)$ and all $f$ in
 ${\cal S } (\R^d)$:
 $$ < Op_h ^{Weyl} (F) f , f> \ \geq - C h^2 \Vert f \Vert _{L^2(\R^d)}^2.
 \leqno (1.4) $$

\bigskip

See [MAR] for these semiclassical versions.
Sometimes the  non negativity of $F$ implies the exact  non negativity of 
the operator, for example in the simple case when $F$ depends on $x$
 or on $\xi$ only. It is possible, too, to apply Melin's inequality. To take 
only one example, let $F\geq 0$ attain its minimum only once,  for a 
nondegenerate critical point. In this case (and in other analogous situations),
Melin's inequality ensures the exact  non negativity of  $Op_h^{Weyl}(F)$ for
a sufficiently small $h$. See  [B-N] or [L-L] for cases when the difference
between  $F(x , \xi)$ and its minimum is equivalent to a power,
 greater than $2$, of the distance between $(x , \xi)$ and the unique
 point where the minimum is attained.

\bigskip

In this article we are interested in the case when $F$ is radial. We assume
 that there exists a function $\Phi$ defined on $\R$ such that
$$ F(x , \xi) = \Phi ( |x|^2+|\xi|^2)   \hskip 2cm (x , \xi) \in \R^{2d}. \leqno (1.5)$$
Moreover, we suppose that $\Phi$ is nondecreasing on  $[0, \infty)$ and 
 such that $F$ is smooth, with bounded derivatives.

In this case, we aim at giving an explicit lower bound 
on the spectrum of the operator $Op_h^{Weyl}(F)$. The main result of
this paper is the following theorem.

\bigskip

{\bf Theorem 1.1} {\it Let $F$  be a smooth function defined on $\R^{2d}$,
 bounded as well as all its derivatives. Assume that $F$ is  of the form 
(1.5),
where  $\Phi$ is a non decreasing function defined on $[0, \infty)$.

Then for all  $f$ in ${\cal S}(\R^d)$,
$$ < Op_h ^{Weyl} (F) f , f> \ \geq
{1\over h} \int_0^{\infty} \Phi (t) e^{-{t\over h}} dt \ \ \Vert f \Vert _{L^2(\R^d)}^2 .
 \leqno (1.6) $$
}

\bigskip
{\bf Remarks}

1 - We do not need  to assume that 
 $\Phi \geq 0$ to ensure the  non negativity of the operator. The
  non negativity
 of the integral suffices.

2 - In the case when $\Phi$ is not flat at the origin, let $m\geq 1$ be the
 smallest integer for which $\Phi^{(m)}(0)\neq 0$. Then one can see that
$$
{1\over h} \int_0^{\infty} \Phi (t) e^{-{t\over h}} dt
= \Phi(0)+ \Phi^{(m)}(0)h^m + {\cal O}(h^{m+1}).
$$

3 - The result can be applied to symbols $F$ depending on the distance from 
another point $(x_0,\xi_0)$ for, if
$\tau F(x,\xi)= F(x+x_0,\xi + \xi_0)$ and $Tf(u)= e^{i (\xi_0 /h)(u-x_0)}f(u-x_0)$, then
$$
 < Op_h ^{Weyl} (\tau F) f , g> \ =  < Op_h ^{Weyl} (F) Tf , Tg> .
$$

\bigskip

We are greatly indebted to N. Lerner for the reference [A-G].

\bigskip

{\bf 2.  Proof of Theorem 1.1.}

\bigskip

We denote by $(H_n)_{(n\geq 0)}$ the sequence of the Hermite functions.
It is a Hermitian basis of
$L^2(\R)$, satisfying
$$ (D^2 + x^2) H_n = (2n+1) H_n. \leqno (2.1)$$

\bigskip
For each multi-index 
 $\alpha = (\alpha _1 , ... \alpha_d)$, we set :
$$ u_{\alpha }(x) = \prod _{j=1}^d H_{\alpha_j }(x_j) . \leqno (2.2) $$
These functions form a Hermitian basis of
 $L^2(\R^d)$.

\bigskip

We shall need the Laguerre polynomials as well, which are defined by
$$ L_n (x) = {e^x \over n!} {d^n \over dx^n} \left ( x^n e^{-x} \right ). \leqno (2.3) $$
One has :
$$ L_0 (x) =1 \hskip 2cm L_1(x) = 1 -x \hskip 2cm L_2(x) =
 {x^2 \over 2 } - 2x +1.   \leqno (2.4) $$

\bigskip

Theorem 1.1 is a consequence of the following proposition, in which 
the parameter $h$ is equal to $1$ and the Weyl operator 
$Op_1^{Weyl} (F)$ is denoted by $Op^{Weyl}(F)$.

\bigskip

{\bf Proposition 2.1} {\it
Under the hypotheses of Theorem 1.1 one has, for all multi-indices 
 $\alpha$ and $\beta $ such that $\alpha \not= \beta$:
$$ < Op^{Weyl}(F) u_{\alpha } ,  u_{\beta }  > =0  .\leqno (2.5) $$
For each multi-index $\alpha $:
$$ < Op^{Weyl}(F) u_{\alpha } ,  u_{\alpha }  > = 2^{-d}  \left [ \Phi (0) V_{\alpha} (0) + {1\over 2} \int_0^{\infty} \Phi'( t/2) V_{\alpha} (t) dt \right ], \leqno (2.6) $$
with 
$$ V_{\alpha }(X) = 4 e^{-{X\over 2}}  \sum _{k=0}^{d-1} C_{d-1}^k  T_{|\alpha |+k} (X), \leqno (2.7)$$
where we set, for all integer  $n$,
$$ T_n (X) = \left [  \sum _{k=0}^{n-1}
 (-1)^k L_k (X) \right ]  + {(-1)^n \over 2} L_n (X). \leqno (2.8) $$

}

\bigskip

{\it Proof of (2.5).}
Let 
 $\alpha $ and $\beta$ be two different multi-indices and let
 $j\leq d$ be such that  $\alpha_j \not= \beta_j$.
Set $P_j = D_j^2 + x_j^2$. According to (2.1) we have :
$$ 2 (\alpha_j - \beta _j) < Op^{Weyl}(F) u_{\alpha } ,  u_{\beta }  > =
< Op^{Weyl}(F) P_j u_{\alpha } ,  u_{\beta }  > -< Op^{Weyl}(F) u_{\alpha } , P_j  u_{\beta }  >.  $$
The fact that $F$ is radial implies that
$x_j {\partial F\over \partial \xi_j}- \xi_j{\partial F\over \partial x_j}= 0$
which, in turn, implies that 
 $Op^{Weyl}(F)$ and $P_j$ commute, thanks to properties of the Weyl calculus.
Consequently, the right term of the above inequality is equal to $0$,
 which proves
 (2.5).

\bigskip

{\it Proof of (2.6).}  
For each multi-index $\alpha$, the Wigner function $ H ( u_{\alpha} , u_{\alpha})$
 (where the parameter $h$, equal to  $1$, is omitted),
satisfies:
$$ H ( u_{\alpha} , u_{\alpha}) (x , \xi) = 2^d (-1)^{|\alpha|}
 e^{-(|x|^2 + |\xi|^2)} \prod _{j=1}^d L_{\alpha _j} (2( x_j^2 + \xi_j^2)). \leqno (2.9)$$
See, for example, [FO] or  [J-L-V]. 
Hence, if $F$ is as in Theorem 1.1,
$$ <Op^{Weyl}(F) u_{\alpha}, u_{\alpha}> =(2\pi )^{-d}  2^d  (-1)^{|\alpha|}
 \int _{\R^{2d} } \Phi (|x|^2+|\xi|^2) e^{-(|x|^2 + |\xi|^2)}
 \prod _{j=1}^d L_{\alpha _j} (2( x_j^2 + \xi_j^2))  dxd\xi. $$
The change of variables $t_j = 2 (x_j^2 + \xi_j ^2)$ allows to write :
$$ <Op^{Weyl}(F) u_{\alpha}, u_{\alpha}> = (2\pi )^{-d}  2^d (-1)^{|\alpha|} (\pi /2)^d
\int _{[0, \infty )^d} \Phi ((t_1 + ... + t_d)/2) e^{-{1\over 2}(t_1 + ... + t_d)}
\prod _{j=1}^d  L_{\alpha _j} ( t_j) dt_1 ... dt_d. $$
This equality can be written as
$$ <Op^{Weyl}(F) u_{\alpha}, u_{\alpha}> =(2\pi )^{-d}  2^d (\pi /2)^d   \int _0^{\infty } \Phi (X/2) U_{\alpha }(X) dX  , $$
with :
$$U_{\alpha }(X) =   (-1)^{|\alpha|} e^{-{X\over 2}}
 \int _{\Omega _d(X)}L_{\alpha _d} (X - t_1 - ... - t_{d-1}) \prod _{j=1}^{d-1}
  L_{\alpha _j} ( t_j) dt_1 ... dt_{d-1}, $$
where
$$ \Omega _d(X) = \{ ( t_1 , ... , t_{d-1} ) , \hskip 1cm t_j >0, \hskip 1cm
t_1 + ... + t_{d-1} < X \}.$$
The equality (2.6) will be a consequence of an integration by parts
 using the following lemma.

\bigskip

{\bf Lemma 2.2} {\it We have:
$$ U_{\alpha }(X) = - V'_{\alpha }(X)  \leqno (2.10)$$
where  $ V_{\alpha }$ is defined by  (2.7) and  (2.8).

}

\bigskip

{\it Proof of Lemma 2.2.} One knows  (cf [M-O-S], section 5.5.2) that
$$ \int _0^X L_{\alpha _1 } (t) L_{\alpha _2 } (X - t)  dt =
L_{\alpha_1 + \alpha_2} (X) - L_{\alpha_1 + \alpha_2+1 } (X).  \leqno (2.11) $$
It follows, by induction on $d$, that
$$ \int_{\Omega _d(X)}  L_{\alpha _d} (X - t_1 - ... - t_{d-1}) \prod _{j=1}^{d-1}
  L_{\alpha _j} ( t_j) dt_1 ... dt_{d-1}
 = \sum _{k=0}^{d-1} C_{d-1}^k (-1)^k L_{|\alpha | + k} (X). $$
Hence
$$U_{\alpha }(X) =     e^{-{X\over 2}}
\sum _{k=0}^{d-1} C_{d-1}^k (-1)^{|\alpha | + k}  L_{|\alpha | + k} (X). $$
Using the recurrence relation $L'_{k+1}(t) = L'_k (t) - L_k (t)$,
we prove (for example by induction)  that for all integer $n$:
$$ {d\over dt } e^{-{t\over 2}} T_n (t) = {(-1)^{n+1} \over 4} L_n (t) e^{-{t\over 2}}
.$$
The equality  (2.10) of the Lemma follows from (2.7) and from the above
 identities.

\bigskip

{\it End of the proof of Theorem 1.1.}
We shall begin by proving (1.6) for $h=1$. Set
$$ S_n (X) = \sum _{k=0}^n (-1)^k L_k(X). \leqno (2.12) $$
Using the recurrence relation
 $L'_{k+1}(t) = L'_k (t) - L_k (t)$, one verifies, by induction,
that for all  $n$:
$$ T'_n (X) = {1\over 2 } S_{n-1}(X). $$
Since $L_n (0) = 1$ for all $n$, we see that $T_n (0) = 1/2$
 and that
$$ T_n (X) = {1 \over 2} + {1 \over 2} \int_0^X S_{n-1} (t)\  dt. $$
According to  [A-G], Theorem 12 (see  [F] as well ),
 $S_n (X) \geq 0$ for all $n\geq 0$ and for all  $X \geq 0$. Therefore
 $T_n (X) \geq 1/2$ for all  $n$ and  $X$,
and, using  (2.7):
$$ V_{\alpha}(X) \geq 2^d e^{-{X \over 2}}. \leqno (2.13)$$
Since $T_n (0) = 1/2$, $ V_{\alpha}(0)  =2^d $.
Hence, if  $\Phi' \geq 0$, one gets :
$$  \Phi (0) V_{\alpha} (0) + {1\over 2} \int_0^{\infty} \Phi'( t/2) V_{\alpha} (t) dt
\geq 2^d \int _0^{\infty} \Phi (t) e^{-t} dt. \leqno (2.14) $$
The inequality
 (1.6), for $h=1$, follows from (2.5), (2.6)  and (2.14). For an arbitrary
 $h>0$, it suffices to apply the above result to the function
 $F_h (x , \xi) = F ( h^{1/2} x , h^{1/2} \xi)$,
that is to say, to the function $\Phi_h (t) = \Phi (th)$.

\bigskip

{\bf References.}

\bigskip

[A-G] R. Askey, G. Gasper, {\it Positive Jacobi polynomial sums, II}, American J. of Math.,
{\bf 98}, 3 (1976), 709-737.

\smallskip

[B-N] R. Brummelhuis,  J. Nourrigat, {\it  A necessary and sufficient condition for Melin's inequality for a class of systems,} J. Anal. Math. 85 (2001), 195–211.

\smallskip

[C-R] M. Combescure, D. Robert, {\it Coherent states and
applications in mathematical physics,} Theoretical and Mathematical
Physics. Springer, Dordrecht, 2012.

\smallskip

[F-P] C. Fefferman,  D. H. Phong, {\it The uncertainty principle and sharp
G\aa rding inequalities, }  Comm. Pure Appl. Math. 34 (1981), no. 3, 285–331.

\smallskip

[F] E. Feldheim, {\it D\'eveloppements en s\'erie de polyn\^omes de Hermite et de
Laguerre \`a l'aide des transformations de Gauss et de Hankel, III,} Konin. Neder. Akad.
van Weten., {\bf 43} (1940), 379-386.

\smallskip

[FO] G. B. Folland, {\it Harmonic Analysis on phase space,} Annals of
Mathematics studies {\bf 122}, Princeton University Press, Princeton
(N.J), 1989.

\smallskip

[HO] L. H\"ormander, {\it The analysis of linear partial
differential operators,} Volume III, Springer, 1985.

 \smallskip

[J-L-V] E.I. Jafarov, S. Lievens, J. Van der Jeugt, {\it The Wigner distribution function
for the one-dimensional parabose oscillator}, J. Phys. A 41 (2008), no. 23, 235301,
and  arXiv:0801.4510.

\smallskip

[L-L] B. Lascar, R. Lascar, {\it In\'egalit\'e  de Melin-H\"ormander en
caract\'eristiques multiples}, to appear, Israel Journal of Math.

\smallskip

[L] N. Lerner, {\it Metrics on the phase space and non-selfadjoint pseudo-differen\-tial operators},
Pseudo-Differen\-tial Operators. Theory and Applications, 3. Birkh\"auser Verlag, Basel, 2010.

\smallskip

 [L-N] N. Lerner, J. Nourrigat, {\it Lower bounds for pseudo-differential operators,} Ann. Inst. Fourier (Grenoble) 40 (1990), no. 3, 657-682.

\smallskip

[M-O-S] W. Magnus,  F. Oberhettinger, R.P.  Soni,
{\it Formulas and theorems for the special functions of mathematical physics.}
Third edition. Die Grundlehren der mathematischen Wissenschaften, Band 52 Springer-Verlag New York, Inc., New York 1966.

\smallskip

[MAR] A. Martinez, {\it  An introduction to semiclassical and microlocal analysis,}
 Universitext. Springer-Verlag, New York, 2002.

\smallskip

[M] A. Melin, {\it Lower bounds for pseudo-differential operators,} Ark. Mat. 9, 117–140. (1971).

\smallskip

[U] A. Unterberger, {\it Les op\'erateurs m\'etadiff\'erentiels},
in Complex analysis, microlocal calculus and relativistic quantum
theory, Lecture Notes in Physics {\bf 126} (1980) 205-241.

\bigskip

Laboratoire de Math\'ematiques, FR CNRS 3399, EA 4535, Universit\'e de Reims
Champagne-Ardenne, Moulin de la Housse, B. P. 1039, F-51687
Reims, France, 

{\it E-mail:} {\tt laurent.amour@univ-reims.fr }

{\it E-mail:} {\tt lisette.jager@univ-reims.fr}

{\it E-mail:} {\tt jean.nourrigat@univ-reims.fr}

\end